\documentclass[twoside,leqno]{article}

\usepackage[letterpaper]{geometry}

\usepackage{siamproceedings,tikz}

\usepackage[T1]{fontenc}
\usepackage{amsfonts,xcolor}
\usepackage{graphicx}
\usepackage{epstopdf}
\usepackage{enumitem}
\usepackage{algorithmic}
\ifpdf
  \DeclareGraphicsExtensions{.eps,.pdf,.png,.jpg}
\else
  \DeclareGraphicsExtensions{.eps}
\fi

\newsiamremark{remark}{Remark}
\newsiamremark{hypothesis}{Hypothesis}
\crefname{hypothesis}{Hypothesis}{Hypotheses}
\newsiamthm{conjecture}{Conjecture}
\newsiamthm{example}{Example}
\newsiamthm{problem}{Problem}
\usepackage{amsopn}

\begin{document}

\newcommand\relatedversion{}
\renewcommand\relatedversion{\thanks{The full version of the paper can be accessed at \protect\url{https://arxiv.org/abs/0000.00000}}} % Replace URL with link to full paper or comment out this line

\title{\Large At the intersection of Numerical Analysis and Spectral Geometry\relatedversion}
    \author{Nilima Nigam\thanks{Department of Mathematics, Simon Fraser University (\email{nigam@math.sfu.ca},  \url{http://www.math.sfu.ca/~nigam}).}
    }

\date{}

\maketitle

% Copyright Statement
% When submitting your final paper to a SIAM proceedings, it is requested that you include
% the appropriate copyright in the footer of the paper.  The copyright added should be
% consistent with the copyright selected on the copyright form submitted with the paper.
% Please note that "20XX" should be changed to the year of the meeting.

% Default Copyright Statement
\fancyfoot[R]{\scriptsize{Copyright \textcopyright\ 20XX by SIAM\\
Unauthorized reproduction of this article is prohibited}}

% Depending on which copyright you agree to when you sign the copyright form, the copyright
% can be changed to one of the following after commenting out the default copyright statement
% above.

%\fancyfoot[R]{\scriptsize{Copyright \textcopyright\ 20XX\\
%Copyright for this paper is retained by authors}}

%\fancyfoot[R]{\scriptsize{Copyright \textcopyright\ 20XX\\
%Copyright retained by principal author's organization}}

%\pagenumbering{arabic}
%\setcounter{page}{1}%Leave this line commented out.

\begin{abstract} How do the geometric properties of a domain impact the spectrum of an operator defined on it? How do we compute accurate and reliable approximations of these spectra? The former question is studied in spectral geometry, and the latter
is a central concern in numerical analysis. In this short expository survey we revisit the process of eigenvalue approximation, from the perspective of computational spectral geometry. 
Over the years a multitude of methods - for discretizing the operator and for the resultant discrete system - have been developed and analyzed in the field of numerical analysis. High-accuracy and provably convergent discretization approaches can be used to examine the interplay between the spectrum of an operator and the geometric properties of the spatial domain or manifold it is defined on. While computations have been used to guide conjectures in spectral geometry, in recent years approximation-theoretic tools and validated computations are also  being used as part of proof strategies in spectral geometry. 
	  Given a particular spectral feature of interest, should we discretize the original problem, or seek a reformulation? Of the many possible approximation strategies, which should we choose? These choices are inextricably linked to the objective: on the one hand, rapid, specialized methods are often ideal for conjecture formulation (prioritizing efficiency and accuracy), whereas schemes with guaranteed, computable error bounds are needed when computation is incorporated into a proof strategy. We also review instances where  the demanding requirements of spectral geometry — the need for rigorous error control or the robust calculation of higher eigenvalues \-- motivate new developments in numerical analysis.
\end{abstract}

\section{Introduction.}

Eigenvalue problems arise in many guises and contexts throughout science and engineering and their mathematical study has occupied mathematicians for several centuries. Indeed, {\it 'it is hard to think of any branch of mathematics where eigenvectors and eigenvalues do not have a major part to play' \cite{Gowers2008}.} Very loosely speaking, if we know all the eigenpairs of an operator, we have access to its properties. The spectral geometer will ask: how do the eigenpairs depend on {\it geometric} properties of the domain on which the operator is defined?  The computational mathematician will ask: how do we {\it compute} these eigenpairs, and if approximation was involved, how accurate is the approximation? In this article, we examine a small subsample of topics which are of mutual interest, and hope this serves as an entry-point into the fascinating world of computational spectral geometry.

{\it Spectral geometry} is the study of the relation between the spectra of operators and the geometric properties of the domains on which these operators are defined (for excellent introductions to various topics, see \cite{ dyatlovbook,levitinbook,colboissteklovsurvey,henrot,girouardreview}. Often the operators of interest are elliptic differential operators on manifolds. The major research developments and questions can (somewhat idiosyncratically) be grouped under three headings:
\begin{itemize}
	\item {\it The forward problem:} Given an  operator ${\mathcal{A}}$ on a domain $\Omega$ (or manifold $(\mathcal{M},g))$, what can one say about its eigenvalues and eigenfunctions? For instance, the beautiful experimental work by Chladni \cite{chladni} is best described as 'On the figures obtained by strewing sand on vibrating surfaces, commonly called acoustic figures', \cite{wheatstone}.  Mathematically: Chladni's work relates the shape of nodal lines of a domain \--- contours where the eigenfunctions are zero \---  to the geometric properties of the given domain \cref{fig:chladni1}. Wheatstone also describes the link between nodal lines and the domain's symmetry and size, \cite{wheatstone}. [Napoleon was fascinated by the experiment, and announced a contest in 1809, with a {\it Prix Extraordinaire} to a person who could explain this striking phenomenon. Sophie Germain was the only competition entrant; after  two attempts, in 1816 she became the first woman to be awarded the Grand Prize in Mathematics of the Paris Academy of Sciences for this work,\cite{germain}]
	\begin{figure}
\centering
		\includegraphics[width=0.3\linewidth]{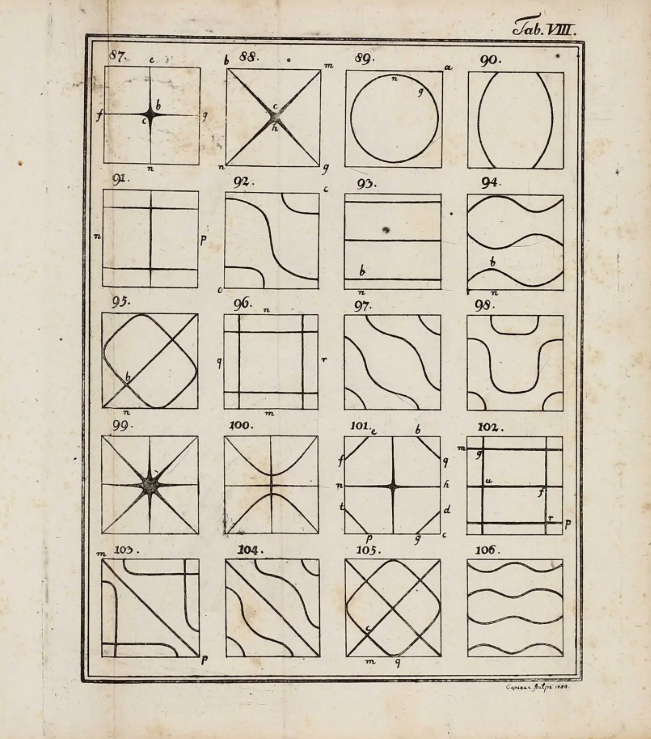}
			\includegraphics[width=0.5\linewidth, height=2in]{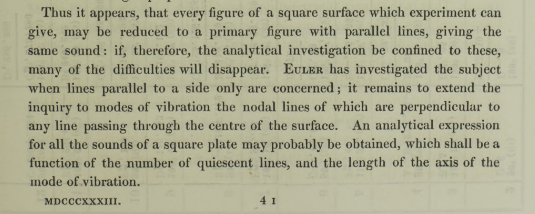}
		\caption{Left: Chladni lines on a square plate from \cite{chladni}. Right: From \cite{wheatstone}}
		\label{fig:chladni1}
	\end{figure} 
	\item {\it The inverse problem:} Given the spectrum of an operator ${\mathcal{A}}$ what can one say about the geometric properties of the domain $\Omega$? The famous question {\it 'Can one hear the shape of a drum?'} due to Kac \cite{Kac1966} falls under this category. Of course, the answer is {\it 'One cannot hear the shape of a drum', \cite{Gordon,Gordon2,buser}.} The existence and construction of such isospectral domains remains an area of active interest. In \cref{fig:isospectral} we show a pair of non-isometric polygons which are provably Dirichlet and Neumann isospectral ( a pair of 'GWW polygons' from the classic construction in \cite{Gordon2,Gordon}). We also show a pair of  planar domains which are isospectral for the mixed Dirichlet-Neumann eigenproblem for the Laplacian, \cite{DNIso}.
	
	\begin{figure}
		\begin{tikzpicture}[scale=2]
			
			% Polygon 1: Vertices (0,0), (1,0), (1.5,0.5), (2,0), (2,1), (1.5,1.5), (0.5,0.5), (0,1)
			\draw[thick, xshift=-1in,fill=gray!30] (0,0) -- (1,0) -- (1.5,0.5) -- (2,0) -- (2,1) -- (1.5,1.5) -- (0.5,0.5) -- (0,1) -- cycle;
			\draw[dashed,xshift=-1in,] (0,0)--(1/2,1/2) --(1,0)--(1,1)--(3/2,1/2)--(2,1)--(1,1);
			% Polygon 2: Vertices (0,0), (0.5,-0.5), (1,0), (0.5,0.5), (1,1), (1,2), (0.5,1.5), (0,2)
			% We use 'xshift=4cm' to move the second polygon to the right.
			\draw[thick,  xshift=0in, fill=gray!50] (0,0) -- (0.5,-0.5) -- (1,0) -- (0.5,0.5) -- (1,1) -- (1,2) -- (0.5,1.5) -- (0,2) -- cycle;
			\draw[dashed,xshift=0in,] (1,0)--(0,0) --(1/2,1/2)--(0,1)--(1,1)--(1/2,3/2)--(0,1);
			%\end{tikzpicture}
			%\begin{tikzpicture}[scale=2.5] % Scale factor for better visibility
			
			% Define the coordinates of the vertices
			\coordinate (A) at (2,0);
			\coordinate (B) at (3.2,0);
			\coordinate (C) at (3.2,1.2);
			\coordinate (D) at (2,1.2);
			
			% Draw the sides with specified colors and a thicker line style
			\draw[line width=2pt,xshift=1in, red] (D)--(A) -- (B) -- (C); % Side BC (Red)
			\draw[line width=2pt,xshift=1in, blue] (C) -- (D); % Side CD (Blue)
			\draw[<->, thin, gray] (2,-0.1) -- (3.2,-0.1) node[midway, below] {$1$};
			\pgfmathsetmacro{\len}{sqrt(2)} % len is approx 1.414
			
			% Define the coordinates (B at origin for right angle)
			\coordinate (B) at (4,0);
			\coordinate (C) at (4+\len,0); % (\sqrt{2}, 0)
			\coordinate (A) at (4,\len); % (0, \sqrt{2})
			
			% Draw the sides with specified colors and a thicker line style
			\draw[line width=2pt, red] (C)--(A) -- (B); % Side AB (Red)
			\draw[line width=2pt, blue] (B) -- (C); % Side BC (Blue)			
			% Draw a square symbol for the right angle at B
			\draw[thick] (0.2, 0) -- (0.2, 0.2) -- (0, 0.2);
			\draw[<->, thin, gray] (4,-0.1) -- (4+\len,-0.1) node[midway, below] {$\sqrt{2}$};
			% Label the vertices
			%\node[below left] at (B) {$B$};
			%\node[below right] at (C) {$C$};
			%\node[above left] at (A) {$A$};
			
		\end{tikzpicture}
		\caption{Isospectral domains. Dirichlet and Neumann isospectral regions (grey polygons)   (Figure 15 \cite{Gordon}). The Dirichlet spectrum of Polygon A  coincides with that of Polygon B.  The Neumann spectra coincide as well. Right: Isospectral domains with mixed Dirichlet (red) and Neumann(blue), Figure 1 from \cite{DNIso}. }\label{fig:isospectral}
	\end{figure}
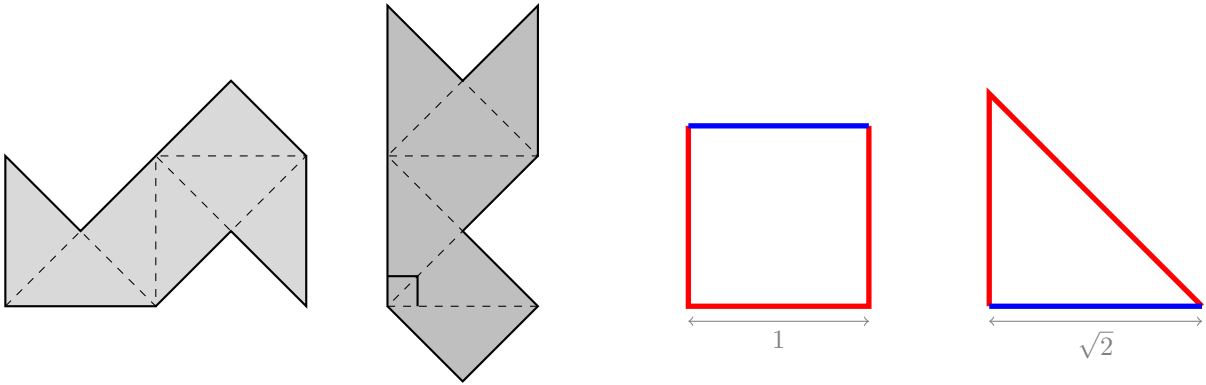
	\begin{figure}
		\centering
		\includegraphics[width=0.15\linewidth]{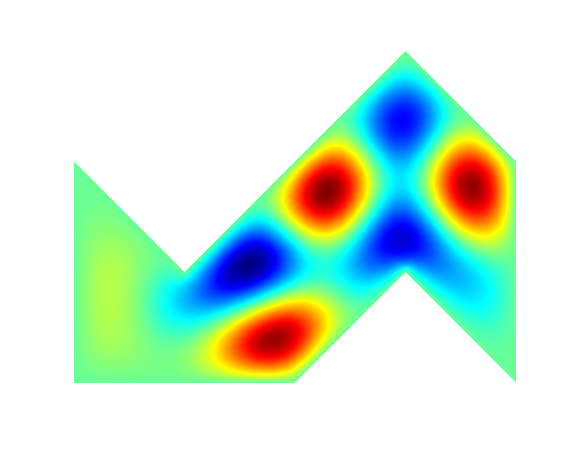}
		\includegraphics[width=0.15\linewidth]{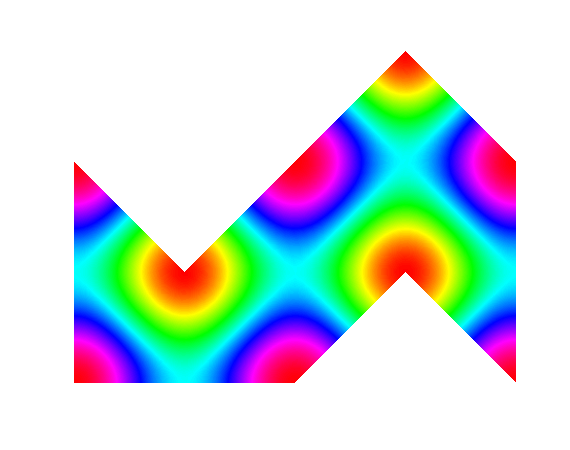}
			\includegraphics[width=0.2\linewidth]{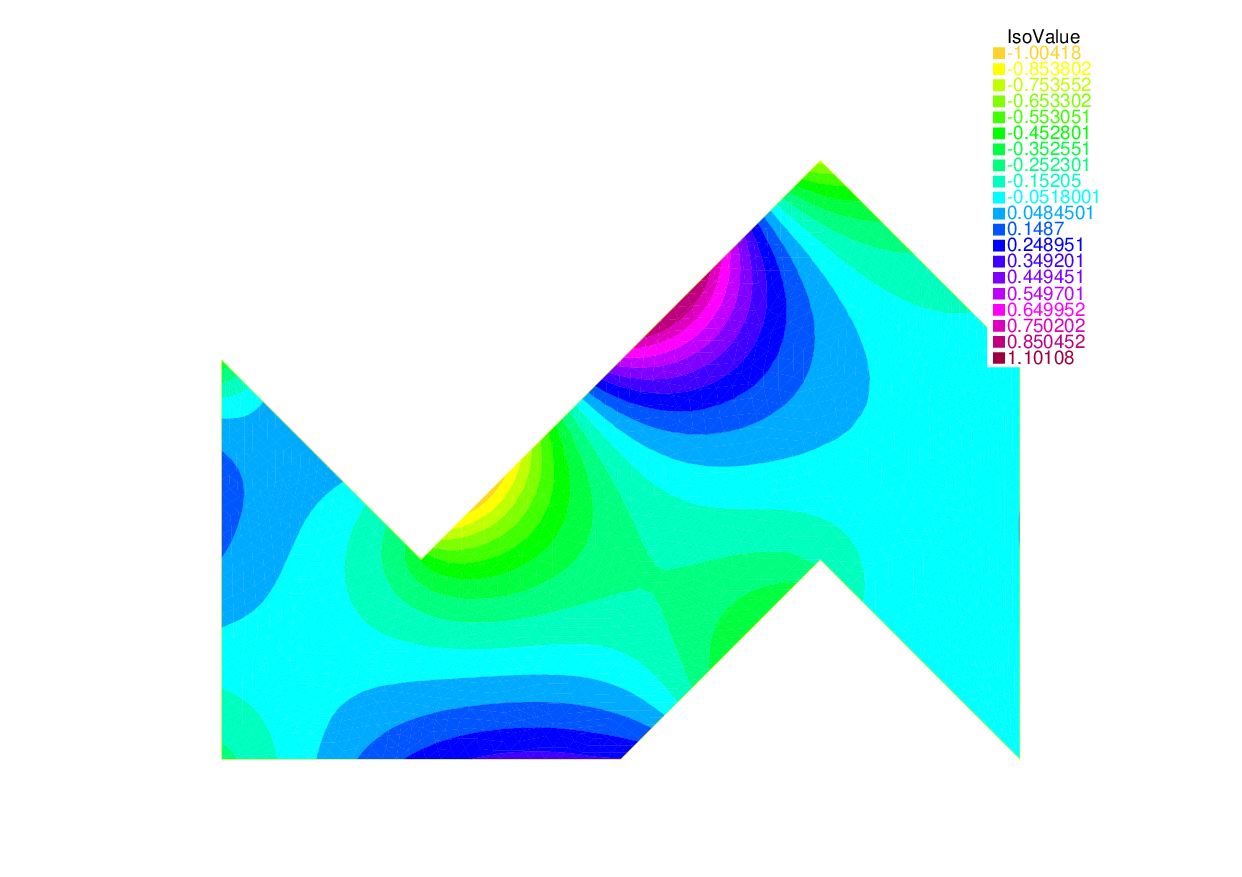}
			
			\includegraphics[width=0.2\linewidth]{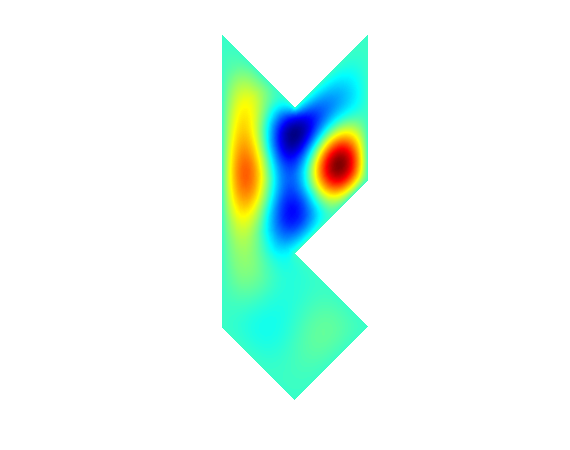}
		\includegraphics[width=0.2\linewidth]{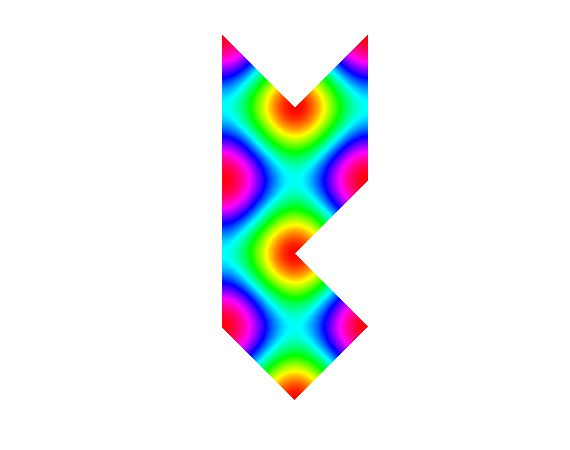}
		\includegraphics[width=0.2\linewidth]{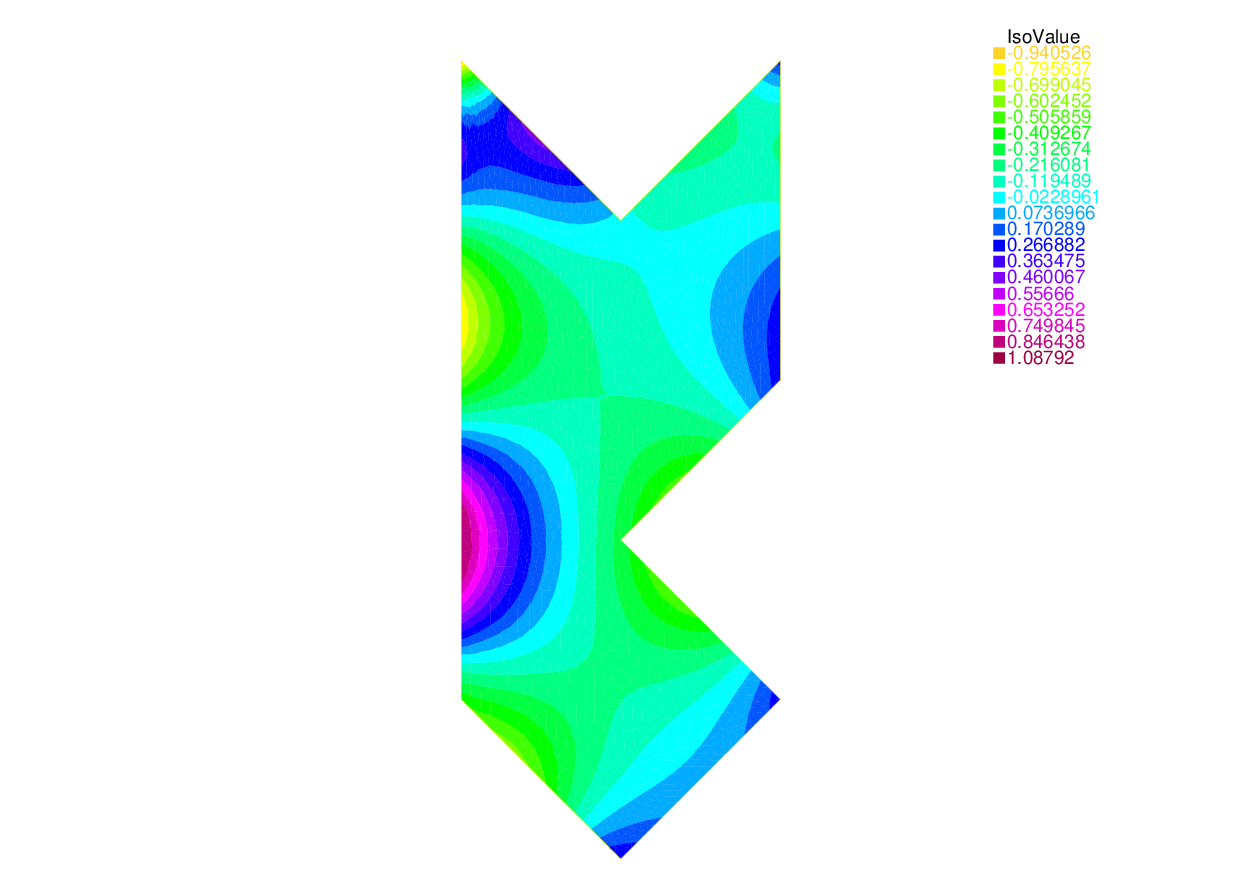}
		\caption{Eigenmodes on the GWW polygons of \cite{Gordon}. Left column, 10th Dirichlet eigenfunction, Dirichlet eigenvalue $\approx 26.08$, computed using the Method of Particular Solutions. Polygons in row 1 and 2 are Dirichlet-isospectral. Middle column: 10th  Neumann eigenfunction, Neumann eigenvalue $\approx 9.62165$, computed using a Boundary Integral approach. Right column: 10th Steklov eigenfunctions, computed using a finite element approach.These domains are not Steklov isospectral, see \cref{table:polysteklov}.  }
		\label{fig:polygonadirichlet9}
	\end{figure}
	
	\item {\it The optimization problem:} Amongst all domains $\Omega \in \mathcal{O}$ of some given class, which one optimizes a given spectral quantity? A classic example is a conjecture from 1877 due to Lord Rayleigh (Chapter IX,\cite{rayleigh}):
		{\it "If the area of a membrane be
		given, there must evidently be some form of boundary
		for which
		the pitch (of the principal tone) is the gravest possible, and this
		form can be no other than the circle."} In other words, the disk minimizes the first Dirichlet eigenvalue amongst bounded planar domains of given area. The statement, independently proven by  Faber \cite{Faber1923} and Krahn \cite{krahn}, can been generalized to bounded domains $\Omega\in \mathbb{R}^n$ and is known  as the Rayleigh-Faber-Krahn inequality:
		\begin{equation}
	\text{argmin}_{\Omega} \left[|\Omega|^{2/n}\lambda_1(\Omega)	\right] = B_1(0,n)\end{equation}
	where $B_1(0,n)$ is the unit ball in $\mathbb{R}^n$.
	%\geq (\text{volume of the unit ball in } \, \mathbb{R}^n)^{2/n}j_{n/2-1,1},
	%	\end{equation} where $j_{m,1}$ is the first non-zero root of $J_m(x)$. 
		The text \cite{henrot} provides an excellent starting point to this field.
\end{itemize}

In numerical analysis, questions about approximation and computation of eigenpairs have been central, for finite and infinite dimensional operators alike. The connection between domain geometry and eigenmodes has been exploited for decades in the sciences and engineering literature and more recently in computer graphics, and numerical approaches are standard here. Precisely because eigenvalue problems are ubiquitous, it is challenging to precisely locate a 'starting point' for their numerical analysis (ad we understand the term today). Suffice it to say: approximating eigenpairs has been a vital part of approximation theory and modern numerical analysis \cite{Courant, collatz}.

More recently, we see the rise in the use of tools from numerical analysis to suggest and prove conjectures in spectral geometry, and this survey attempts to review some of these developments. The numerical analyst has a wide array of discretization approaches to choose from, each with advantages and disadvantages. Should one use a variational method, or a collocation approach? Are piecewise polynomial approximation spaces preferred, or are global bases better? Should one always use gradient-based methods for optimization? What is the best approach for solving the discrete system?  

 These questions take on a new light when used to formulate conjectures, or as part of proofs. Which of a plethora of problem formulations and numerical methods is suitable for a given question in spectral geometry? Which questions in spectral geometry are amenable to computation? How does the goal - formulating conjectures or proving statements - affect these choices? In this paper, we aim to provide an entry point to these questions.  
 
 The intersection of numerical analysis and spectral geometry includes many fascinating problems, and perspectives from both fields are needed.  An example of the interplay between these fields is \cite{becker}, where the  spectral characterization of 'magic angles' in twisted bilayer graphene is studied. The mathematical model leads to non-Hermitian operators whose spectrum is both rich and difficult to simulate; novel numerical methods which carefully account for pseudospectral effects are devised. In another example, a novel and accurate computational strategy for the Neumann-Poincar\'e operator provided strong numerical evidence for the existence of embedded eigenvalues \cite{helsing1} within a continuous spectrum (the domain is an ellipse with a corner attached); the existence of such eigenvalues was  proven in \cite{shipman}. 
 
 As these examples highlight, the world of  spectral geometry problems of interest to numerical analysts is considerably bigger than the small corner we examine here. 

\section{A bird's eye view of spectral computations.}

We begin with an abstract discussion of the process of computing approximate eigenpairs. While the process is standard, it is nonetheless useful to examine the steps involved in discretization \--- as we shall see, different choices will dictate if the computed approximations are sufficiently accurate as to be used to help formulate conjecture or provide evidence. Other choices may be necessary if we use the approximation procedure as part of a proof.

	To set notation, consider the following spectral problem: 
	\begin{problem}\label{problem:true} Find 
		$	\tilde{u} \in \tilde{\mathcal{H}}, \lambda \in \mathbb{C}$ so that
		\begin{equation}\label{eq:true}
			\tilde{\mathcal{A}} \tilde{u}= \lambda\tilde{u}.
	\end{equation}\end{problem}Here $\tilde{\mathcal{H}}$ is a suitable (infinite dimensional) space of functions on $\Omega$, and $\tilde{\mathcal{A}}:\text{dom}(\tilde{\mathcal{A}})\subset\tilde{\mathcal{H}}\rightarrow \tilde{\mathcal{H}}$.

At this juncture, it is a good idea to ask if this problem is amenable to solution via numerical approaches.
\subsection{Can we always compute the spectrum of an operator?}
If by 'compute' we mean the use of approximations, the answer is 'no'. It is not at all obvious that the discrete spectrum computed by a seemingly-reasonable discretization approach preserves the multiplicity of the eigenspace, or indeed that 'converged' approximate eigenvalues are actually in the true spectrum. The question of computability is even more challenging if the spectrum is continuous or essential.  'Thar be monsters' is an apt description; certainly the 'traps and snares' of approximation are well-known within the numerical analysis literature (see, e.g.\cite{boyd}). These issues can be particularly severe while discretizing nonlinear spectral problems, and this remains a very interesting and active area of research. However, even for linear operators, it is easy to find examples where 'reasonable' truncation or discretization approaches fail.  We provide some cautionary examples below.

\begin{example} \cite{davies,hansen,bottcher} Let $T_\epsilon:\ell^2(\mathbb{Z}) \rightarrow \ell^2(\mathbb{Z})$ be defined by a variant of the shift operator,
	$$ T_\epsilon(f)(n):=\begin{cases}
		\epsilon(n+1),\qquad n=0\\
		f(n+1).\qquad \text{otherwise}
	\end{cases}$$
We easily confirm for $\epsilon>0$,  the spectrum $\text{spec}(T_\epsilon)=\{w\in \mathbb{C}\vert |w|=1\}$ while for  $\epsilon=0$, $\text{spec}(T_\epsilon)=\{w\in \mathbb{C}\vert |w|\leq 1\}$. A 'small perturbation' of the operator completely changes the spectrum! 
	\end{example}
This illustrative example of spectral instability, attributed to E.B. Davies \cite{davies} (older references exist, see \cite{bottcher}) raised the fundamental question: {\it Is it even possible to {\it compute} the spectrum and pseudospectrum of an arbitrary closed linear operator on a separable Hilbert space?} Suppose one cannot distinguish a small number $\epsilon>0$ from 0(as may be the case while working in finite precision). The spectrum of $T_0$ cannot be approximated via that of $T_\epsilon$. Does this render the project of spectral approximation hopeless? This long-standing question is answered (thankfully we can compute the spectrum in many cases!) in \cite{hansen}, but as the paper highlights, these issues are both challenging and subtle. [In the case of this particular example, one instead approximates the {\it n-pseudospectrum:} a set which has the continuity properties of the $\epsilon-$pseudospectrum of a closed and bounded operator, but which also approximates the original spectrum,\cite{hansen}.] Even simple questions \--- for which classes of differential operators on unbounded domains can one find computational algorithms
that converge to the true spectrum, and also guarantee that the algorithm output is in the spectrum up
to an arbitrary small distance?\--- require careful examination, and were substantially open until recently. This exciting area of computational mathematics and numerical analysis has helped put the questions of spectral approximation on a firmer footing through the introduction of a Solvability Complexity Index,  \cite{hansen,colbrook21,colbrooke23,colbrook24}.  Table 1 in \cite{colbrooke23} summarizes the  situation of spectral approximation through this lens. Of particular relevance in spectral geometry is the recent work \cite{colbrook24} about the computability of geometric features (e.g. Lebesgue measure, spectral radius, capacity) of spectra.

Many fascinating problems in spectral geometry arise in situations where the spectrum includes not just eigenvalues, and these can easily fall under the more challenging categories of the Solvability Complexity Index. At least for this numerical analyst, the possibility of resolving a continuous spectrum or locating embedded eigenvalues is daunting. Even more daunting is the following situation: suppose we compute approximations of the spectrum of an operator $\tilde{\mathcal{A}}$ on a varying domain $\Omega$. Can we tell \--- just from the computed spectrum \--- that the operator has   continuous spectrum (or worse)?  Computational approaches for such problems are deeply challenging, and a subject of intense current research, (see e.g.\cite{helsing1}). 
For the purpose of exposition, we mostly restrict ourselves here to questions concerning the spectra of familiar differential operators on bounded domains in $\mathbb{R}^2$ and $\mathbb{R}^3$, where the spectrum has only discrete and countable eigenvalues. Even for these (simpler) problems, discretization of the operator must be done with care. 

\subsection{Which problem should we solve?}
\label{sec:hyperef}
Suppose the operator $\tilde{\mathcal{A}} =\tilde{\mathcal{A}}(\Omega)$ in \cref{problem:true} is defined on a domain $\Omega\in \mathbb{R}^n$, and has a countable discrete point spectrum. If we want to study the dependance of the {\it eigenvalues} of $\tilde{\mathcal{A}}(\Omega)$ on $\Omega$, we could instead study a related problem whose spectrum is the same, albeit easier/more accurate/faster to compute.
There are many possible reformulations one could consider. In our notation, the reformulated problem is described as follows.
\begin{problem}Find \label{problem:new}
	$ u_i \in \mathcal{H}, \lambda_i\in \mathbb{C}, i=1,2,...$ so that
	\begin{equation}
		\label{eq:new}\mathcal{A}(\lambda_i) u_i =0, 
	\end{equation}
	where $\mathcal{H}$ is a Hilbert space of functions on $\Omega$, and $\mathcal{A}:\text{dom}(\mathcal{A})\subset \mathcal{H}\rightarrow \mathcal{H}.$
\end{problem} 
Let us illustrate this with a simple example.
\begin{example}\label{example:D} (Dirichlet problem for $-\Delta$)
	Let $\Omega\in \mathbb{R}^2$ be a bounded domain with Lipschitz boundary $\partial \Omega$, and define $\tilde{\mathcal{H}}:= C^2(\Omega)\cap C(\bar{\Omega}).$ The Dirichlet eigenvalue problem is: find 
	$	\tilde{u}_i\not=0 \in \tilde{\mathcal{H}} , \lambda_i \in \mathbb{R}, i=1,2...$ so that
	\begin{equation}
		-\Delta \tilde{u}_i= \lambda_i \tilde{u}_i, \qquad \tilde{u}_i\vert_{\partial \Omega}=0, \quad i=1,2,...
	\end{equation}
	Let $\mathcal{H}=$ the Sobolev space $H_0^1(\Omega)$ and $\mathcal{T}_1:\mathcal{H}\rightarrow \mathcal{H}$ be the resolvent operator for the Dirichlet Laplacian with data $f\in H^1_0(\Omega)$ defined by
	$$ \int_\Omega \nabla(\mathcal{T}_1 f) \cdot \nabla \phi = \int_\Omega f \phi, \qquad \forall f \in \mathcal{H}.$$ The operator $\mathcal{T}_1$ is compact and self-adjoint. A reformulation of  the Dirichlet eigenvalue problem is: find $v_i\in \mathcal{H}, \mu_i \in \mathbb{R}, i=1,2....$ so that
	$$\mathcal{A}_1(\mu_i)v_i := (\mathcal{T}_1-\mu_i \mathcal{I})v_i=0.$$
	These eigenvalues are related to the Dirichlet eigenvalues as $\mu_i = \frac{1}{\lambda_i},$
	and the eigenfunctions of $\mathcal{T}_1,\tilde{\mathcal{A}}$ are the same.
	\end{example}
	We note the eigenfunctions $v_i\in \mathcal{H}$ in this reformulation are defined on $\Omega$. Possible reformulations  in terms of functions defined on $\partial \Omega$ are derived using integral operators.
	
	\begin{example} Consider again the Dirichlet problem for the Laplacian on $\Omega$. Let $\mathcal{H}:=H^{-1/2}(\partial \Omega)$. Let $k\not=0 \in \mathbb{R}$, and define the boundary integral operators $\mathcal{S}_k: \mathcal{H} \to H^{1/2}(\partial\Omega)$ be defined as
	\begin{align}
		\mathcal{S}_ku[ \psi](x)  &:= \int _{\partial \Omega}\Phi_k(x,y)\psi (y)ds_y \qquad \mbox{for } x \in \partial \Omega, 
	%	\mathcal{K}_k[ \psi](x)& :=  -\frac{\psi (x)}{2}+\int _{\partial \Omega}\frac{\partial}{\partial
	%		n_x}\Phi_k(x,y)\psi (y)ds_y  \qquad \mbox{for } x \in \partial \Omega.
	\end{align} where $\Phi_k(x,y):=H_0^1(k |x-y|),$ $H_0^1$ the Hankel function of order 1 for wavenumber $k$. Let $\phi_i\not=0 \in \mathcal{H}, k_i\in \mathbb{R}$ solve
	\begin{equation}
		\mathcal{S}_{k_i}[\phi_i](x)=0,\qquad x\in \partial \Omega.
	\end{equation} 
	Then $\lambda_i = k_i^2, i=1,2...$,\cite{AKHMETGALIYEV20151,ammari20}. 
\end{example}
We note the reformulated problem need not be a PDE eigenvalue problem.  Instead, we seek the real characteristic values of the operator $k \rightarrow\mathcal{S}_{k}.$ The Dirichlet eigenfunctions $u_i$ of the Dirichlet problem coincide with the single layer potential corresponding to $\phi_i$, that is, 
\begin{equation} u_i(x)=\tilde{S}_{k_i}[\phi_i](x):= \int _{\partial \Omega}\Phi_k(x,y)\psi (y)ds_y \qquad \mbox{for } x \in \Omega.\end{equation}
Other reformulations in terms of integral operators can be found, and the unknown $\phi_i$ is defined on the boundary. This reduction in domain can prove advantageous.
As a final example involving a reformulation, we consider the Steklov problem for the Laplacian. We note that this could also be reformulated in terms of the Dirichlet-to-Neumann operator.
\begin{example}\label{example:S}(Steklov problem for $-\Delta$.)
Let $\Omega\in \mathbb{R}^2$ be a bounded domain with polygonal boundary $\partial \Omega$, and $\tilde{\mathcal{H}}:=C^2(\Omega)\cap C(\bar{\Omega})$. The Steklov eigenvalue problem is: find 
$	\tilde{u}_i \in \tilde{\mathcal{H}}, \sigma_i \in \mathbb{R}, i=1,2...$ so that
\begin{equation}
	-\Delta \tilde{u}_i= 0, \qquad \frac{\partial}{\partial n}\tilde{u}_i\vert_{\partial \Omega}=\sigma_i \tilde{u}_i\vert_{\partial \Omega}, \quad i=1,2,...
\end{equation}Let $\mathcal{H}:=H^1(\Omega)$, and define $\mathcal{A}_2:\mathcal{H}\rightarrow \mathcal{H}$ as 
$$ \int_\Omega \nabla(\mathcal{A}_2 f) \cdot \nabla \phi + \int_{\partial\Omega} \mathcal{A}_2f \phi = \int_{\partial\Omega} f \phi, \qquad \forall \phi \in \mathcal{H}.$$  $\mathcal{A}_2$ is a compact and self-adjoint operator, \cite{mora}. The spectrum of ${\mathcal{A}}_2 ={0,1}\cup\{\mu_j\}_{j \in \mathbb{N}}$ where $\mu=1$ is an eigenvalue corresponding to the constant function, $\mu=0$ is has an infinite-multiplicity eigenspace $\{w \in \mathcal{H}\vert \int_{\partial \Omega} w=0\}$, and the other eigenvalues have finite multiplicity. The Steklov eigenvalues $\sigma = \frac{1}{\mu}-$ for the non-zero eigenvalues $\mu$.

We can also use an indirect approach via a single layer ansatz to derive an isospectral problem. Define the boundary operators (see \cite{kressbook})
$\mathcal{S}_0: H^{-1/2}(\partial \Omega) \to H^{1/2}(\partial\Omega)$, $\mathcal{K}':H^{-1/2}(\partial \Omega)\rightarrow H^{-1/2}(\Omega)$ as
\begin{align}
	\mathcal{S}_0[ \psi](x)  &:= \int _{\partial \Omega}\Phi_0(x,y)(\psi (y)-\bar{\psi})ds_y \qquad \mbox{for } x \in \partial \Omega, \\
	\mathcal{K}[\psi](x)& :=  \int _{\partial \Omega}\frac{\partial}{\partial		n_x}\Phi_0(x,y)(\psi(y)) ds_y\qquad \mbox{for } x \in \partial \Omega
\end{align} with $\Phi_0(x,y):=\log(|x-y|), \bar{\psi} := \frac{1}{|\partial \Omega|}\int_{\partial \Omega} \psi(y)\, ds_y$. Using standard jump relations (\cite{hsiaobook}), we obtain  the eigenvalue problem: find $\phi_i\not=0\in H^{-1/2}(\Omega),\sigma_i\in\mathbb{R}$ such that
\begin{equation}\label{eq:bieformulate}
\left(\frac{1}{2}\mathcal{I}+ \mathcal{K'}\right)[\phi_i-\bar{\phi_i}] = \mu_i \mathcal{S}_0[\phi_i]\end{equation}
The eigenvalues $\mu_i$ coincide with the Steklov eigenvalues. The  eigenfunctions are related as
$$ u_i(x) = \int _{\partial \Omega}\Phi_0(x,y)(\psi (y)-\bar{\psi})ds_y \qquad \mbox{for } x \in  \Omega.$$ 
\end{example}

As we see, many reformulations of the original spectral problem are possible. It may be possible to identify another, isospectral problem. Or one may be able to map the domain $\Omega$ to another, more computationally tractable one (via conformal mapping, stereographic projection, etc.) Whether to use a reformulation should be guided by the problem, as well as the goal. In the examples above, the operator $\mathcal{A} = -\Delta $ has constant coefficients; the integral reformulation is easily achieved since the fundamental solution is known; the associated operators are very well studied, as is their use for exterior problems (\cite{hsiaobook,Saranen, coltonandkress}). The associated eigenfunction $\phi$ is only defined on $\partial \Omega$, leading to a reduction of dimension from $\mathbb{R}^2$. If our goal is to approximate the eigenvalues of $\mathcal{A}$, this is an excellent approach.  If the operator $\mathcal{A}$ includes a variable coefficient or if the fundamental solution is not available, using integral operators may not be an easy option. 

\subsection{Enter: Approximation.}

We note that until this stage, no approximations have been made. If we seek  approximations to the solutions of \cref{problem:true} or \cref{problem:new} on a computer in finite time, then we cannot expect to locate {\it all} the eigenvalues. We instead seek a finite number of pairs $u_{i,N} \in H_N, \lambda_{i,N} \in \mathbb{C}$, for 
 $N \in \mathbb{N}$ and an $N-$dimensional function space $H_N$. If $H_N\subset \mathcal{H}$ the method is {\it conforming}, otherwise not. The spaces $H_N$ and $\mathcal{H}$ are related through a projection $P_N:\mathcal{H}\rightarrow H_N$, and the specific projection specifies whether the given approach is a Galerkin or collocation approach.  The choice of a particular function space $H_N$ defines the method, for example finite difference methods, the method of particular solutions, or spectral methods. 
\begin{problem}Find \label{problem:approx}
	$ u_{i,N}\in H_N, \lambda_{i,N}\in \mathbb{C}, i=1,2,...,N$ so that
	\begin{equation}
		\label{eq:new}A_N (\lambda_{i,N})(u_{i,N}) = 0,
	\end{equation}
	where $H_N$ is now a finite-dimensional space of functions on $\Omega$, and $A_N:H_N\rightarrow H_N$. 
\end{problem} 
Before we proceed, we ask: in what sense do solutions of \cref{problem:approx} approximate those of \cref{problem:true} or \cref{problem:new}? As we have seen, this question requires careful examination \--- and is much more involved than the analogous question for the source problem $A_N (w)=f$. Indeed, even a discretely stable discretization for the source problem may fail to correctly approximate the eigenspaces. Even worse, in some situations, the discrete eigenvalues converge to {\it spurious} eigenvalues. Without information about the spectral problem and the approximation method, one cannot infer much from the enumerate numerical eigenvalues \--- these may appear to converge, and yet have nothing to do with the original spectrum of interest. We point to the excellent works \cite{boffidefinition,boffi} for deeper discussions and cautionary examples about these issues.

 Suppose we want to approximate $\lambda_k,$ the $kth$ eigenvalue of $\tilde{\mathcal{A}}$, and suppose $\mathcal{A}_N(\lambda_{k,N})u_{k,N} := \mathcal{B}_Nu_{k,N} - \lambda_{k,N}u_{k,N}=0$ Clearly  $\lambda_{k,N}$ is the $kth$ discrete eigenvalue of $\mathcal{B}$. Clearly we need $\lambda_{k,N} \rightarrow \lambda_k$ as $N$ becomes large. We need to ensure $\lambda_{k,N}$ does not converge to a spurious eigenvalue for any $k$ being computed.  The convergence of eigenfunctions must be handled appropriately: if $\lambda_k$ is an eigenvalue of multiplicity $m$, there is no guarantee $\lambda_{k,N}$ {\it also} has the same multiplicity. These considerations lead to the following definition \cite{boffidefinition}:
 
 \begin{definition}
 	Let positive integer $k$ be any positive integer. Let $m(k)$
 	denote the dimension of the space spanned by the $first$ distinct $k$ eigenspaces corresponding to the first $k$ distinct eigenvalues. That is, if $E_i$ denotes the eigenspace associated to $\lambda_i$, define  
 	$$ m(1)=\text{dim}\{\oplus_i E_i, \lambda_i=\lambda_1\},\qquad m(\ell+1) =m(\ell) + \text{dim}\{\oplus_i E_i, \lambda_i=\lambda_{m(\ell)+1}\}.$$ The first $N$ {\it distinct} eigenvalues are $\lambda_{m(1)},\lambda_{m(2)},...\lambda_{m(N)}$.
 	We say the discrete eigenvalue problem \cref{problem:approx} converges to the
 	continuous one \cref{problem:new} if, for any $\epsilon$ and $k>0$, there exists $N_\epsilon\in \mathbb{N}$ such
 	that, for all $N>N_\epsilon$, we have
 	\begin{enumerate}
 		\item {\it Convergence of eigenvalues:}
 	$$ \max_{1\leq i \leq m(k)}|\lambda_i - \lambda_{i,N}|<\epsilon.$$
 	\item {\it Convergence of eigenspaces:}
 	$$ \hat{\delta}\left(\oplus_{i=1}^{m(k)}E_i, \oplus_{i=1}^{m(k)}E_{i}^N\right)<\epsilon.$$
 	 \end{enumerate} Here, the gap $\hat{\delta}$ between two subspaces $U,V$ of a given Hilbert space $H$  is given by 
 	 $$ \hat{\delta}(U,V):=\max \left(\delta(U,V),\delta(U,V)\right),\qquad \delta(U,V):=\sup_{u\in U,\|u\|_H=1}\inf_{v\in V}\|u-v\|_H.$$
 \end{definition}
While this definition is well-known in the numerical analysis literature, we emphasize its importance to problems in spectral geometry: these discretizations will be used either to guide the formulation of conjectures, or as part of proofs. It is critical, therefore, that we design eigenvalue approximation approaches which are convergent according to the definition above.

Assume that \cref{problem:approx} converges to \cref{problem:new}. We still need to specify in what sense \cref{problem:approx} is being solved: we want to minimize the residual $R(w,\ell):=A_Nw - \ell w$. Since\cref{problem:approx} is posed in a finite-dimensional space $H_N$, we can select a basis $\{\phi_i\}_{i=1}^N$ for the trial space $H_N$, and write the eigenfunctions in \cref{problem:approx} as $u_{i,N}:=\sum_{j=1}^Nc^N_{i,j} \phi_j$ for $i=1,2,...k$. Using the {\it method of weighted residuals} we obtain:
\begin{problem}\label{problem:discrete}
Find  $u_{i,N}:=\sum_{j=1}^Nc^N_{i,j} \phi_j \in H_N, \lambda_i\in \mathbb{R}$ so that
\begin{equation}
	| \left\langle \mathcal{A}_Nu_{i,N},\psi - \lambda_{i}^Nu_{i,N},\psi\right\rangle_w| =0 \qquad \forall \psi \in G_N.
\end{equation}
 Here $G_N:=\text{span}\{\psi_j\}_{j=1}^N$ is a test space, and $\langle \cdot,\cdot \rangle_w$ is a suitable bilinear form or duality pairing.
\end{problem}
How we minimize the residual completes the description of our method. We note that the test and trial spaces can be chosen to have different dimension, in which case the system obtained will be rectangular.

\begin{example}\label{example:Dgalerkin} Continuing \cref{example:D}, we seek the eigenvalues of the compact solution operator $\mathcal{A}_1:\mathcal{H}\rightarrow \mathcal{H}.$ Let $H_N\subset \mathcal{H}=H^1_0(\Omega)$ denote a $N-$dimensional conforming space of piecewise $p$-degree polynomial finite elements, corresponding to a triangulation of $\Omega$, and let $W_N = \text{span}\{\psi_i\}_{i=1}^N \subset H^1_0(\Omega)$ consist of some other approximation space.  The method of weighted residuals leads to a Petrov-Galerkin approach to the Dirichlet eigenvalue problem:
	\begin{align}
		\left\langle \mathcal{A}_1 u_i,\psi\right\rangle_w - \lambda\left\langle  u_i,\psi\right\rangle_w   =0 \Leftrightarrow \int_{\Omega} \nabla(\mathcal{A}_1u_i)\cdot \psi_j = \mu_i\int_{\Omega} \nabla u_i \cdot \nabla \psi_j \Leftrightarrow \int_{\Omega} u_i\cdot \psi_j = \mu_i\int_{\Omega} \nabla u_i \cdot \nabla \psi_j,\forall j=1,2,...N.
	\end{align}
	
\end{example} If we select $W_N=V_N$, the discrete problem is symmetric for this example. We could also choose different test and trial spaces \--- as in a Petrov-Galerkin \--- scheme.
We could also devise a pseudospectral/collocation approach for the eigenvalue problem.

\begin{example} In \cref{example:D}, suppose $\Omega$ is a polygon with $M$ corners with interior angles $\pi/\alpha_i, i=1,2,..M$. 
The Fourier-Bessel function $ \phi(r,\theta) = J_{\alpha_i k}(\sqrt{\lambda}r)\sin(\alpha_i k \theta)$ satisfies $-\Delta w=\lambda w $ in $\Omega$, and its trace vanishes on the edges of the polygon forming the $ith$ corner. We form $H_N$ from a collection of such functions. The residual is now the mismatch at the boundary, and is minimized pointwise at boundary collocation points $x_i,i=1..N$:
\begin{align}
	\left\langle u_i,\delta_{x_i}\right\rangle_w -0   \Leftrightarrow \sum_{i=1}^N c_i \phi_i(x_\ell)  = 0,\forall \ell=1,2,...N.
\end{align}
 This is the Method of Particular Solutions, \cite{henrici,betcke}.
\end{example}

Finite element approaches by design are closest to the variational characterization of eigenvalues, and can easily be formulated for complex geometries and variable-coefficient problems. Their convergence behaviour on a range of eigenvalue problems is very well understood (see, eg.,the texts \cite{babuska,boffi,sunbook}). The finite element systems arising in \cref{problem:discrete} typically lead to sparse matrices. If the space $H_N$ consists of piecewise polynomials (as is typical) the rate of convergence depends on the polynomial degree and the regularity of the eigenspaces. The integrals which arise are easily computed. Typical mesh triangulations are via simplices, leading to another layer of approximation on  smooth domains. 

Integral reformulations for elliptic eigenvalue problems have been explored extensively both in the mathematics and the engineering literature, for the Laplacian but also for elastic vibrations, transmission eigenvalue problems, mixed Dirichlet-Neumann eigenvalue problems, Steklov eigenvalues, resonance problems. We point to  \cite{Kitahara,AMINI1995208,lu91,duranBEM,Steinbach, unger14,ammaribook,Cakoni02012017,Akhmetgaliyev2017,AKHMETGALIYEV20151} amongst others. As mentioned, we need a fundamental solution for these reformulations, and the associated eigenfunctions are defined on the boundary. Curvilinear boundaries are easily handled via parametrization, as are mixed boundary conditions or exterior problems. We note that using integral formulations we can use {\it either} a Galerkin approach (BEM) or a collocation approach. 
The associated system in \cref{problem:discrete} leads to smaller but denser matrices than FEM; for smooth boundaries, a collocation approach can provide very high accuracy eigenvalues. Closely related, the Method of Particular Solutions relies on having simple closed-form solutions of the differential operator. They also lead to small and dense matrices, and high-accuracy approximation. As with integral equation methods, these are less adaptable to variable coefficient problems.

Which method should we pick for a given spectral geometry question? The answer depends on the goal: if we seek fast and accurate approximations of eigenvalues for constant-coefficient problems, MPS and integral-operator based methods are excellent. For variable-coefficient or eigenvalue problems with low regularity, the FEM approach is well-suited. If the goal is to use these computations as part of a proof strategy, the finite element approach may offer advantages.

\subsection{Approximation of the spectra of discrete systems}

Having chosen the formulation, the finite-dimensional subspace and the method by which the residual is to be minimized, one arrives at a finite-dimensional system: either generalized eigenvalue problems for the eigenpairs ${\tt \lambda_{i,N},u_{i,N}}$ or the problem of locating ${\tt \lambda_{i,N}}$ so that the parameter-dependent $\tt{A}(\ell)$ is singular when $\ell=\lambda_{i,N}$.  That is, the problems \cref{problem:discrete} are now of the form: 
\begin{problem}\label{problem:gevp}{\it find $\mathtt{\lambda_{i,N}}\in \mathbb{C}$ and nonzero vectors ${\tt u_{i,N}}$ for $i=1,2,...N$ such that } $$ \mathtt{A_N u_{i,N} = \lambda_{i,N} B_N u_{i,N} } \quad \text{or} \quad  \mathtt{A_N(\lambda_{i,N})u_{i,N} =0}.$$ 
\end{problem}
In most instances these problems {\it cannot} be solved exactly, and another layer of approximation occurs. Perhaps the first approximation approach to computing eigenvalues is from 1846, when Jacobi proposes a method to diagonalize a real symmetric matrix \cite{jacobi}. Typically the generalized eigenvalue problems in \cref{problem:gevp} are large, and one selects amongst iterative eigenvalue solvers (such as the Lanczos or Arnoldi methods, or shift-and-invert approaches) best suited to the structure of the generalized eigenvalue problem. In case one seeks the characteristic values of ${\tt A}(\ell),$ it is necessary to devise a fast search method over parameter values $\ell \in \mathbb{C}$; these are typically based on contour integral approaches (See, for instance, \cite{sunbook,jayjeff}). 

The accurate approximation of the solutions of \cref{problem:gevp} is deeply challenging, particularly in the presence of clusters or multiple eigenvalues. The literature on the approximation theory of discrete eigenvalue and resonance problems is truly vast, and we do not touch upon it, pointing instead to some excellent expository references \cite{arnoldi,GOLUB200035, saad,embreebook}.

\section{Computational Spectral Optimization}
In recent years numerical approaches for spectral shape optimization in $\mathbb{R}^d$  have played a particularly important role in formulating conjectures. 
We point to the recent papers
\cite{Ammari,ammari20,oudet21,MR3640624,MR4052734,MR4058408,Osting2010,Osting2012,Osting2017a,Viator2022, kao2025extremalsteklovneumanneigenvalues} amongst others, for a glimpse into this exciting area of research in computational spectral geometry. As an example, in \cite{MR3640624} consider the problem of maximizing the $k$ Steklov eigenvalue amongst fixed-volume sets in $\mathbb{R}^d$. In addition to theoretical results, the paper includes some fascinating conjectures, based on numerical experiments. One such: amongst all fixed-area planar domains in $\mathbb{R}^2$, the maximizers of the $kth$ eigenvalue are connected, and have $k-fold$ symmetry. These computations used the MPS to compute eigenvalues, and the results agreed with those using a boundary integral approach, \cite{Akhmetgaliyev2017}.

The typical computational spectral shape optimization question for a spectral quantity $\mathcal{J}_k(\tilde{\mathcal A},\Omega)$  involving the atmost the first  $k$ eigenvalues of $\tilde{A}(\Omega)$ can phrased as: amongst all domains $\Omega$ in some admissable class $\mathcal{D}$, find 
\begin{equation} \text{arg opt}_{\Omega\in \mathcal{D}} \mathcal{J}_k(\tilde{\mathcal A},\Omega)\end{equation}
There are three main challenges for these problems: the parametrization of domains $\Omega\in \mathcal{D}$ using atmost a finite number of parameters; the construction of a shape derivative to locate descent directions and how to handle topological changes. 

During the optimization process, one must compute the (discrete approximations of) the eigenvalues of $\tilde{A}(\Omega)$ on the candidate domain, and the considerations discussed in the previous section become even more relevant, with speed and efficiency becoming an added issue. One may opt for a blend of methods: a fast eigenvalue method is used for the initial part of the optimization, and then as one approaches the candidate optimizer, one can use more accurate discretizations.

Working with these computational constraints can lead to new mathematical results in spectral geometry. In \cite{ammari20}  we studied a spectral optimization problem for a mixed Dirichlet-Neumann eigenvalue problem. The optimization was over length and location of the Neumann portions of the boundary with the constraint that one of the eigenvalues was close to a prescribed target value. As part of the strategy, we derived asymptotic expressions for the perturbation of the mixed Dirichlet-Neumann eigenvalues, when a small portion of the boundary is changed from Dirichlet to Neumann. This allowed us to reduce the total number of eigenvalue solves during the optimization procedure. We studied a similar optimization problem for the mixed Steklov-Neumann problem in \cite{Ammari}, where again the location and length of the Neumann piece is found.  In very recent work,  \cite{kao2025extremalsteklovneumanneigenvalues} study this maximization problem over partitions of fixed measure. 
%As a final comment: for many spectral optimization problems even the existence of a global optimizer is hard to establish; the role of numerical experiments in this case is to identify possible local optimizers. 
\section{Spectral Inverse Problems }
The area of inverse problems in spectral geometry is truly vast, and we do not attempt a summary. Instead, we focus on Kac's question: can we always distinguish shapes given their spectra? The question of isospectrality provides a nice illustration of the interplay between the fields of spectral geometry and numerical analysis.  Some of the known Dirichlet isospectral domains are planar, non-convex polygons \--- and are {\it proven} to be isospectral \cite{Gordon,Gordon2}.  For numerical analysts, these examples are irresistible: these are easy-to-construct domains  therefore serve as tests of accuracy of numerical eigenvalue approaches. As a very nice example, in \cite{DriscollTobinA.1997EoID} the author significantly improves the accuracy of a previous method \cite{Descloux} which was based a variational approach via Fourier-Bessel functions.  The same isospectral domains could also be studied via the collocation-based MPS based on the Fourier-Bessel functions(method originally due to \cite{henrici}). MPS and potentially able to provide eigenvalue bounds \cite{henrici} \--- but is known to be highly ill-conditioned in the presence of re-entrant corners such as the isospectral drums.  The accurate MPS computation of the spectra of these drums is revisited in a another famous paper \cite{betcke}, in which the authors propose an elegant fix to the ill-conditioning.  

In the other direction, what can computations tell us about isospectrality?  If the discrete spectra agree up to a finite number of eigenvalues would provide evidence (but not a proof) of isospectrality. It is easier to show two domains are not isospectral.
\begin{example} Do the GWW polygons have the same Steklov spectrum? In \ref{table:polysteklov} we present numerical results via several different methods. These results provide {\it evidence} but not {\it a proof} that the GWW polygons are not Steklov-isospectral. We thank David Sher and Carolyn Gordon for suggesting this example.
\begin{table}\centering
	\begin{tabular}{|c|c|c|c|c|c|}\hline
		&k&P1nc&BIO&P2&P1\\ \hline
PolyA&1&.2748&    0.2795&    0.2803&    0.2845\\
PolyB& & 0.3033 &  0.3088 &   0.3096 &   0.3144\\ \hline
PolyA&2&0.7804&    0.7905 &   0.7919 &   0.8014\\
PolyB&&0.6015 &   0.6117  &  0.6130  &  0.6244 \\ \hline
PolyA& 3&1.0812  &  1.0892  &  1.0897  &  1.0980\\
PolyB & & 1.2147  &  1.2364  &  1.2375  &  1.2624\\ \hline
PolyA& 4&1.6751&    1.7030  &  1.7054 &   1.7331	\\
PolyB&& 1.9995  &  2.0239 &   2.0244 &   2.0477 \\ \hline
	\end{tabular}\label{table:polysteklov}
	\caption{First 4 non-zero Steklov eigenvalues for the GWW polygons of \cref{fig:isospectral}. These polygons are clearly {\it not} Steklov isospectral. We show only the first 5 digits. The computations are via finite element approaches: P1 conforming method (P1), a P2 conforming approach (P2) and the Crouzeix-Raviart elements (P1nc). Computations via a boundary integral collocation approach (BIO) are also shown. Note the P1 and P1nc values appear to bracket the other (more accurate) ones. }
\end{table}
\end{example}

\section{Guaranteed eigenvalue bounds}
A central and recurring preoccupation in numerical analysis is: for a given discretization, how close to the true spectrum (for us, a discrete and countable set) is the (finite) set of computed eigenvalues? As seen in previous sections, one works with approximation strategies for which the discrete spectrum is convergent, in the sense that $\|\lambda_{k,n} - \lambda_k|\rightarrow 0$ as the discretization parameter $n\rightarrow \infty$. This is an asymptotic result. If computations are to be used in a proof, we must address the question: for fixed $N$, how close is $\lambda_{k,N}$ to $\lambda_k$?

In a seminal paper by 
\cite{henrici}, a much-used eigenvalue bound for the Dirichlet eigenvalues of a linear self-adjoint second-order elliptic operator $D$ on  $\Omega \subset\mathbb{R}^n$ is established using a maximum principle argument, which is particularly useful when the eigenpairs are approximated via the Method of Particular Solutions or via layer potentials.
\begin{theorem}\cite{henrici} Suppose
  $(u_{h},\lambda_{h})$ is a computed eigenpair of operatpr $D$, with $\|u_h\|_{L^2(\Omega)}=1$.  Further assume $Du_{h} = \lambda_{h} u_{h}$ in $\Omega, $ with $\max_{x\in \partial \Omega}|u_{h}(x)|<\epsilon<1$.  Then there is a true Dirichlet eigenvalue $\lambda$ of $D$ on $\Omega$ 
  such that
  $$ |\lambda-\lambda_h| \leq \lambda_h \left(\frac{\sqrt{2}\epsilon+\epsilon^2}{1-\epsilon^2}\right).$$
\end{theorem}
It is a nontrivial task to estimate $\max_{x\in \partial \Omega}|u_{h}(x)|,$ and in the original paper by Fox, Henrici and Moler, the authors estimate this term for $-\Delta$ on a planar polygonal domain by means of a line search over each edge. 
Eigenvalue bounds in a similar spirit for the Neumann problem were recently derived in \cite{NeumannBounds}; these in turn allow for accurate calculations of eigenvalues far into the spectrum.

It is easy to see that conforming finite element approaches for eigenvalues will yield approximations from above, ie., $\lambda\leq \lambda_{k,h}$; bounds of the form $|\lambda_k -\lambda_{k,h}| \leq C_{k}(\Omega) h^{r}$ can be derived where $h$ is the mesh size and $r$ depends on the polynomial degree and the regularity of the eigenspace. Unfortunately, the constant $C_k(\Omega)>0$ in these theorems is typically not computable. Explicit {\it guaranteed lower bounds} via $\lambda_{k,h}$, if available, can provide a powerful strategy for spectral geometry questions involving eigenvalue inequalities.

 As a concrete instance, in  \cite{jakobson} it was conjectured that the first eigenvalue of the Laplacian on a surface of genus two was maximized by a singular metric. Using symmetries, the eigenvalue problem on the surface is reduced to a certain mixed Dirichlet-Neumann boundary value problem on a half-disk:
\begin{align}\label{eq:dima}
	-\Delta u = \lambda \frac{4}{(1+r^2)^2} u, \, x\in \Omega, &\qquad  u \vert_{\Gamma_D}=0, \, \frac{\partial u}{\partial n}\vert_{\Gamma_N} = 0.
\end{align}
It remained to show that the first eigenvalue $\lambda>2.$ We used both conforming and non-conforming finite elements to approximate this eigenvalue and {\it observed} numerically that using conforming FEM, $\lambda_{h} \searrow 2.27...$ whereas with non-conforming FEM $\lambda_{h}\nearrow 2.27...$. While the numerical evidence was compelling, at the time there were few available results in the numerical analysis literature on guaranteed lower bounds for eigenvalues. One notable exception was \cite{armentano}, in which the following theorem was shown:
\begin{theorem} (Theorem 2.3 of \cite{armentano}).\label{armentano}
	Let $\lambda_k$ be the $k$th Laplace-Dirichlet eigenvalue on a planar polygonal domain $D$, and let $\lambda_{k,h}$ be an approximation computed using a nonconforming Crouzeix-Raviart space with zero Dirichlet condition, $V_h$ ($h>0$ the mesh size). If the Dirichlet eigenfunction $u_k$ belongs to the Besov space  $B_2^{1+r,\infty}(D)$ and if there is a constant $c>0$ such that 
	$\|u_k- u_{k,h}\|_h\geq c h^r, r<1$ then for $h$ small enough
	$$ \lambda_{k,h}<\lambda_k.$$
\end{theorem}
Unfortunately, this theorem was not directly applicable to the situation of \cref{eq:dima}. result is asymptotic in $h$, and the constant $c$ is not explicitly computable. Moreover, the domain $\Omega$ under consideration is not polygonal. As a result, the statement $\lambda>2$ for the eigenvalue problem \cref{eq:dima}  remained a conjecture until it was proven 12 years later (\cite{NAYATANI201984}, preprint 2017). 
\begin{figure}
	\centering
	\includegraphics[width=0.15\linewidth]{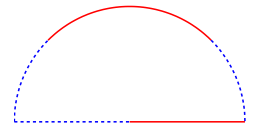}
	\caption{Half-disk $\Omega$ for problem \cref{eq:dima}. Dirichlet conditions are imposed on the red curves, and Neumann conditions on the blue curves. From \cite{jakobson}.}
	\label{fig:dima}
\end{figure}

In parallel developments, other  research groups were independently pursuing the question of {\it computable} lower bounds for eigenvalues of the Laplacian on polygonal domains using finite element discretizations. The papers \cite{oishi13} and \cite{carstensen} were the first to provide computable lower bounds, with \cite{liubook} focussing on conforming approximations, and \cite{carstensen} on the use of non-conforming methods. In both cases, the lower bounds involved are computable. As a prototypical result, assuming the discrete eigenvalue problem using Crouzeix-Raviart elements are exactly solved, we obtain 
\begin{theorem} (Theorem 3.2 of \cite{carstensen}). Let $\Omega$ be a planar polygonal domain.
Let $(\lambda_h,u_j) \in \mathbb{R}\times V_h$ 
be the non-conforming approximation of the Dirichlet eigenpair $(\lambda,u)$ of the smallest eigenvalue. Then 
$ \frac{\lambda_{h}}{1+\kappa^2h^2\lambda_h}\leq \lambda,$ with $\kappa^2 := (\frac{1}{8} + \frac{1}{j_{1,1}^2})$ and $h$ the maximal mesh size. 
\end{theorem} 
The proof of this theorem relies in part on a  estimate of the $L^2$ interpolation error using the Crouzeix-Raviart elements on triangular elements. This in turn relies on a Poincar\'e inequality on triangles \cite{laugesen} \--- which is obtained by solving a spectral optimization problem on triangles! Theorem 3.1 in the same paper provides a computable lower bound which includes the impact of approximating the eigenpairs of the discrete generalized eigenvalue problem. Papers in this line of work include \cite{xie2018,carstensen25} and references therein.

Beginning with \cite{oishi13,Liu2015}, other groups have derived guaranteed eigenvalue bounds via conforming methods, and used them in combination with validated numerics and interval arithmetic to prove conjectures in spectral geometry.  We point to the very recent paper \cite{endo} concerning the simplicity of the second Dirichlet eigenvalue for nearly degenerate triangles. This paper highlights the rich cross-fertilization between the fields of numerical analysis and spectral geometry. In \cite{Liu2015} computable lower bounds of Laplacian eigenvalues are evaluated by means of a projection
error estimate for a nonconforming finite element space; in \cite{endo}, this approach is generalized to infinite dimensions, to provide a computer-assisted proof of

\begin{theorem} (Theorem 1.1 of \cite{endo}) The second Dirichlet eigenvalue is simple for every non-equilateral triangle with its minimum
	normalized height less than or equal to $\frac{1}{2}\tan(\pi/60)$.
\end{theorem} 
 \section{Numerical analysis, spectral geometry and conjectures}
 The use of careful computations to provide insight and help formulate conjectures is a not-surprising form of interaction between the fields of numerical analysis and spectral geometry. Provided the numerical strategy used is provably stable and convergent, we can use the computed eigenpairs \--- even after only a finite number of mesh refinements or iterations \--- to guide our mathematical ideas.
 \subsection{How numerical experiments can lead to a conjecture}s
  We consider the Steklov spectrum on an annular domain $\Omega_\epsilon$ given by the region between the unit circle centered at the origin $\partial B_1(0,0)$, and the circle $\partial B_{0.1}(0,\epsilon)$ with $ \epsilon\in [0,0.88])$. This is a problem with a countable, discrete spectrum.  Since we are interested in the eigenvalues, we consider instead the reformulation in  \cref{eq:bieformulate} (after suitably modifying the jumps). This leads to a discrete problem where the unknowns are only on the boundary. The problem is discretized using a well-known collocation approach, \cite{coltonandkress}. The discrete eigenvalue problem is solved using the QZ algorithm.
	
\begin{figure}
		\centering
		\includegraphics[width=0.4\linewidth]{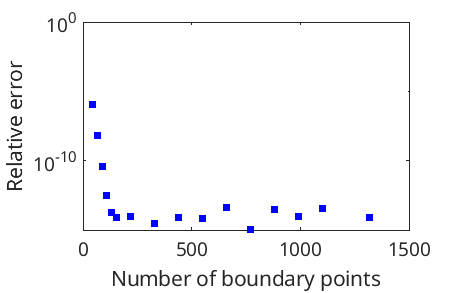}
		\includegraphics[width=0.4\linewidth]{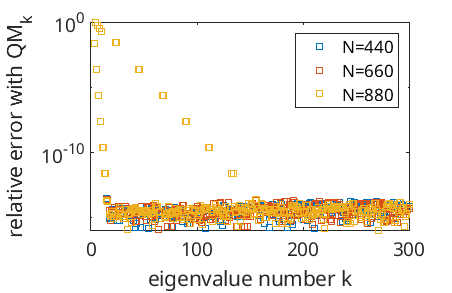}
		\caption{Steklov eigenvalues on an annulus $\Omega_0$. Left: the relative error of the first 20 nonzero eigenvalues, $\sum_{n=1}^{20} \frac{\sigma_k^{true} - \sigma^N_k}{\sigma_k^{true}}$, as a function of the number of discretization points on the boundary. Right: the relative error $\frac{\sigma_k^{N} - S_k}{S_k}$ for eigenvalue numbers $k=3...400$, with $N=440, 660$ or $880$ collocation points.   }
		\label{fig:concentricconvergencerelerror}
	\end{figure}

	We note that for $\Omega_0$ (formed by concentric circles), one can find the Steklov eigenvalues in closed form via separation of variables (the resultant expressions involve powers of $\epsilon^n$). This allows us to check the accuracy of our numerical approach, see the left subplot of \cref{fig:concentricconvergencerelerror}. The eigenvalues are seen to converge spectrally. We highlight the fact that $\sigma_0=0$ and $\sigma_{9} =-(1.1)/(0.1log(0.1)) \approx4.777239300935770$ are eigenvalues of multiplicity 1, and the rest are double eigenvalues.  
	
\begin{figure}[H]
		\centering
		\includegraphics[width=0.4\linewidth]{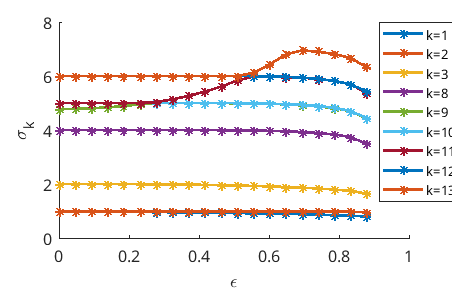}
		\includegraphics[width=0.4\linewidth]{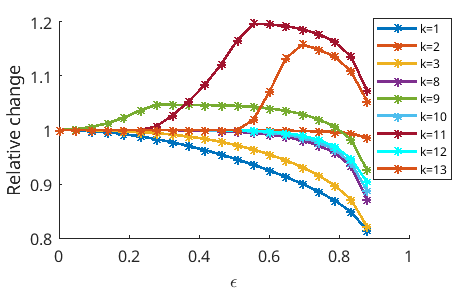}
		\caption{How does the Steklov spectrum of an annulus vary as the inner disk moves? The center of the inner disk in $\Omega_\epsilon$ is at $(0,\epsilon)$. Left: $\sigma_k(\Omega_\epsilon)$ v/s $\epsilon$ for $k1,2,3,8,9,10,11,12,13$. Right: relative change $\sigma_k(\Omega_\epsilon)/\sigma_k(\Omega_0)$.}
		\label{fig:sigmarandom}
	\end{figure}

	Following the results of \cite{alexspectral}, we expect the Steklov spectrum $ \{\sigma_k\}$ of $\Omega_\epsilon$ to converge rapidly to that of the union $S$ of the Steklov spectra of $B_1(0,0)$ and $B_{0.1}(0,0)$. If we denote $S:=\{S_k\}_{k=1}^\infty$ the latter set (reordered),  $\sigma_k(\Omega_{\epsilon}) = S_k + \mathcal{O}(k^{-\infty}).$ This provides a qualitative check of the computed eigenvalues: we expect the computed eigenvalues $\sigma_k^N$ to agree with $S_k$ for large $k$. In \cref{fig:concentricconvergencerelerror} we observe $\sigma^N_k \rightarrow S_k$ as $k$ increases, including two subsequences which converge slower than the others.	
	We now ask: how does $\sigma_k(\Omega_{\epsilon})$ depend on the location of the inner circle as $\epsilon$ varies in $[0,0.88)$? In \cref{table:annulus} we first document the performance of the discretization approach for $\epsilon=0.88$; in this case, the distance between the inner and outer boundary is $0.02$. As can be seen, the method retains high accuracy with a modest number of discretization points, even for the larger eigenvalues. In Figure \cref{fig:sigmarandom}, we see how $\sigma_k(\Omega_\epsilon)$ varies as the position of the inner disk changes.
	
	\begin{table}[H]\centering
 	{\small 	\begin{tabular}{|c|c|c|c|c|c|} \hline
			N&$\sigma_1$& $\sigma_2$&$\sigma_{10}$ & $\sigma_{100}$\\ \hline
			130&0.792792051250894 &   0.961874727024116 &   4.443052025754944 & 46.610918905024320\\
			260&0.794544789783840 &   0.961790500646763 &   4.438623167632640 &46.449879036723786\\
			520&0.794597544967050 &   0.961791478912838 	&  4.438646392806247 & 46.438542901400155\\
			780&0.794597555469579 &   0.961791479149678  & 4.438646399420676& 46.438543189261750\\
			1040&0.794597555472255&  0.961791479149744  & 4.438646399422233&  46.438543189337942 \\ \hline 
		\end{tabular}	\label{table:annulus}}
		\caption{Computed Steklov eigenvalues on $\Omega_{\epsilon}, \epsilon=0.88$ for $k1,2,10,100$. Here, $N$ is the number of quadrature points on the boundary. }
	\end{table}	
 The accuracy of our approach has been checked against true eigenvalues where they are known. In the case the true eigenvalues are not known, the numerical convergence of the computed values is checked. Another form of validation (via the 'quasimodes' of \cite{alexspectral}) is also provided in an 'extreme' situation. We are confident the numerical approach yields accurate results, and allow us to formulate a conjecture.  
 
\begin{conjecture}
	Let $B_r(a,b)$ denote the ball of radius $r>0$ centered at $(a,b)$, and denote
 $\Omega_e:=B_{1}(0,0)\setminus B_{\delta}(0,\epsilon)$. The first non-zero Steklov eigenvalue on $\sigma_1(\Omega_{\epsilon})$ is a monotone decreasing function of $\epsilon, \epsilon\in [0,1-\delta)$.
\end{conjecture}
	We thank Alexandre Girouard for bringing this example to our attention.

\subsection{Numerical approaches and conjectures in spectral geometry}
This is a non-exhaustive list of examples where careful numerical approximations (via discretization of PDE or integral operators) are used to formulate conjectures in spectral geometry. We focus mainly on situations combining finite element/boundary integral techniques.

\begin{itemize}
	\item FEM approaches (\cite{oudet}) as well as MPS approaches (\cite{antunesfreitas}) have been used to conjecture minimizers of $kth$ Dirichlet eigenvalues and maximizers of Neumann eigenvalues for domains with fixed area, $k\geq 3$, including evidence of optimizing non-symmetric shapes eg $\lambda_7$ for the Dirichlet problem.
	\item Numerical evidence via boundary integral approaches supports  conjectures regarding fixed-area domains maximizing the $kth$ Steklov eigenvalue for the Laplacian in the plane \cite{Akhmetgaliyev2017}, including the uniqueness and $k-fold$ symmetry of such optimizers. In $\mathbb{R}^3$ and $\mathbb{R}^4$, MFS approaches supports  conjectures regarding fixed-volume domains maximizing the $kth$ Steklov eigenvalue for the Laplacian in $\mathbb{R}^3,\mathbb{R}^4$, \cite{antunessteklov}. In some instances the numerical optimizers do not have $k-fold$ symmetry.
	\item FEM are used for a relaxation of the maximization problem for the $kth$ Neumann eigenvalue problem, leading to candidate domains\cite{bucurdensity}.
	
	\item Conjecture supported by FEM calculations regarding the optimizing  metric $g$ for the first Laplace-Beltrami eigenvalue on a closed surface of genus 2 with area fixed \cite{jakobson} (proved in \cite{NAYATANI201984}).
%	\item Polya's conjecture for the Aharonov–Bohm operator on unit disk \cite{Frank}. FEM approaches were used. Proved in \cite{filinovbohm}.
	\item Numerical evidence in support of a Poly\'a-Szeg\H{o} conjecture for Steklov on $n-$gons, \cite{sebasteklov}, using FEM combined with Bayesian optimization. 
		
	\item FEM computations were used to conjecture  formulae for the scattering matrix for some hyperbolic surfaces of genus one \cite{strohmeiercusp}.
	
	\item MPS-based computations of the eigenvalues of  the Aharanov-Bohm on disks,  annuli and square with centered potential were presented in (\cite{Frank}); this lead to a Poly\'a-style inequalities which were later proved - with a computer-assisted  component - in (,\cite{filinovballs,filinovbohm},\cite{filonov2025polyasconjecturedirichleteigenvalues}).
	\item Numerical observations generated via the MFS regarding Steklov-Lam\'e maximizers in the plane (volume constraints), \cite{ANTUNES20251}, including: the maximizer of $\lambda_1(\Omega)$ appears to be  the disk.
	\item Numerical experiments via boundary integral approaches \cite{oscarjeff} answer (in the negative) 
	{\it Open Problem 10 \cite{girouardreview}. (i) Are the nodal sets of Steklov eigenfunctions on a Riemannian manifold $\Omega$ dense on the scale $\frac{1}{\sigma}$ in $\Omega$?}. In the same paper, this result is proved.	
\end{itemize}

\section{Numerical analysis, spectral geometry and proof}
When discretization approaches are used as part of a mathematical proof strategy \--- even for a small part \--- we must pay attention to all steps of the discretization procedure.  We may need to 'prepare' the proof in surprising ways. As an example, in \cite{nigamfocm} we prove a modification of Schiffer's conjecture on a regular pentagon. As part of the proof, we need to use finite element computations on certain subdomains of a triangle $T$ with vertices $ (0,0),(1,0)$ and $(0,\tan(\frac{\pi}{5}))$. Since one vertex of this triangle is not rational, we needed to introduce  other triangles, whose vertex coordinates {\it are} rational, allowing for mass and stiffness matrices with rational entries. To finish the proof guaranteed eigenvalue bounds were used. However, in order to work on these triangles instead of the original $T$, we needed to prove related domain monotonicity results. See  \cite{mesh} for similar discussions.

\subsection{Computations and proof: some examples}	

This is a non-exhaustive list of examples where techniques from numerical analysis and approximations (via discretization of PDE or integral operators) are used to prove conjectures in spectral geometry. We focus mainly on situations combining finite element, MPS or boundary integral techniques with rigorous computing. There are other examples of rigorous computing being used to prove theorems in spectral geometry (see, e.g, \cite{bourque2023linearprogrammingboundshyperbolic}) which we do not discuss here.

\begin{itemize}
	\item The second Dirichlet eigenvalue on nearly-degenerate triangles is simple, \cite{endo}. The key conjecture (Conjecture 6.47 of \cite{henrot}) is proved using rigorous numerics, expanding on ideas on FEM-based eigenvalue bounds.
	
	\item Polya's conjecture for the Dirichlet problem on triangles, using guaranteed FEM bounds \cite{endo23} 
	\item The Neumann-Polya conjecture for eigenvalues $\lambda \in [3,14]$ of the unit ball in $\mathbb{R}^2$ Theorem 4.4 of \cite{filinovballs} proved via rigorous numerics.
	
	\item  Polya conjecture for Dirichlet eigenvalues of the Laplacian on annuli proven in \cite{filonov2025polyasconjecturedirichleteigenvalues}; Theorem 8.1 proved via rigorous numerics. 
	\item A modification of Schiffer's conjecture on regular pentagons was proved via a finite element approach on rational meshes, \cite{nigamfocm}
	
	\item In 1997, \cite{levitinnodal} demonstrated a counterexample of Payne's nodal line conjecture (that the nodal line of  the second Dirichlet eigenfunction on $\Omega$  must touch the boundary of $\Omega$ ), by means of a planar, bounded, non-simply connected domain for which the nodal line is closed and does not touch the boundary. In turn, they asked: {\it "What is the smallest $N$
		such that there exists a domain with $N$ 
		boundary components whose second eigenfunction has a nodal line that does not hit the boundary"?} This question is partially  addressed by Theorem 1.2 in \cite{DAHNE2021105957}: {\it There exists a planar domain with 6 holes for which the nodal line of $u_2$ 
		is closed;} the theorem is proved using by MPS and FEM, coupled with validated numerics.
		\item  The result {\it There exist two triangles $A$ and $B$, not isometric to each other, such that $\lambda_i(A) =
		\lambda_i(B),  i = 1, 2, 4$} is proved in  \cite{GOMEZSERRANO2021920} using  a combination of FEM+ MPS + rigorous numerics are used as part of the proof strategy.
	
	\item In \cite{bogoselbucur1} the authors provide a proof strategy for the The Polya-Szeg\"o conjecture on polygons: {\it 'Amongst n-sided planar polygons of fixed area, the first Dirichlet eigenvalue is minimized by the regular n-gon'}. The proof uses atmost a finite number of numerically-certified computations. They then prove the result \--- assuming exact arithmetic and exact solution of the discrete problems arising from finite element computations \--- for $n=5,6,7,8$. Validated numerics with interval arithmetic are used in \cite{bogosel2024polygonalfaberkrahninequalitylocal} to demonstrate the local optimality of the regular pentagon and hexagon.
	
\end{itemize}
\section{Conclusion}
Spectral approximation can be fruitfully used both to provide insight and as part of proof strategies for problems in spectral geometry. These efforts in turn  help drive developments in the numerical analysis of spectral problems. 

 Amongst the possibilities for future collaboration between spectral geometry and numerical analysis  we mention the spectral geometry of vectorial problems on domains (and generalizations to forms/manifolds), the continued development of computable eigenvalue and eigenspace bounds particularly in the presence of clusters or multiple eigenvalues, for domains with curvilinear boundary or on smanifolds, metric optimization. Guaranteed methods in spectral optimization to certify optimizers and inverse problems as well as in the study of operators with complicated spectra would also valuable both mathematically and in terms of algorithmic progress. Structure-preserving discretizations, probabilistic approaches and validated numerics are likely to play an ever-increasing role as well. We foresee many continued collaborative efforts on challenging problems at the interface of numerical analysis and spectral geometry. 

\section*{Acknowledgments.}
The author is grateful to the Canadian Natural Sciences and Engineering Research Council of Canada. Thanks to Eldar Akhmetgaliyev, Oscar Bruno, Habib Ammari, Nurbek Tazhimbetov and Kshitij Patil for many fruitful discussions and code concerning  boundary integral eigenvalue problems. The author is grateful to Dorin Bucur, Michael Levitin, Iosif Polterovitch and David Sher. Finally, Alexandre Girouard and Tom Archibald are sincerely thanked fors helpful discussions about this article. 

% SIAM recommends using BibTeX
% if using BibTeX
\bibliographystyle{siamplain}
\bibliography{Nigam_ICM}
\end{document}